\documentclass[12pt]{article}
\usepackage[latin1]{inputenc}
\usepackage[english]{babel}
\usepackage{fancybox}
\usepackage{amssymb, amsmath} 
\usepackage[a4paper, centering]{geometry}
\usepackage{graphicx,color}
\usepackage{subfigure}
\usepackage[T1]{fontenc}
\usepackage{amsthm}
\makeatletter
\def\th@plain{%
  \thm@notefont{}
  \itshape 
}
\def\th@definition{%
  \thm@notefont{}
  \normalfont 
}
\makeatother
\newcommand{\res}{\mathop{\hbox{\vrule height 7pt width .5pt depth 0pt
\vrule height .5pt width 6pt depth 0pt}}\nolimits}

\newtheorem{thm}{Theorem}[section]
\newtheorem{prop}[thm]{Proposition}
\newtheorem{lem}[thm]{Lemma}

\theoremstyle{definition}
\newtheorem{defn}[thm]{Definition}
\newtheorem{example}[thm]{Example}
\newtheorem{oss}[thm]{Remark}

\def\proof{\noindent {\bf Proof.}\quad\ignorespaces}
\def\endproof{{\hfill\qed}\par}

\newcommand{\Z}{\mathbb{Z}}
\newcommand{\N}{\mathbb{N}}
\newcommand{\R}{\mathbb{R}}
\renewcommand{\epsilon}{\varepsilon}

\DeclareMathOperator{\dist}{dist} 

\begin{document}

 \vspace{1.5cm}
                   \begin{center}
   {\large \bf VARIATIONAL PROBLEMS WITH PERCOLATION:\\[0.1cm]
                  RIGID SPIN SYSTEMS}\\[1cm]
                   \end{center}
\begin{center}      
                  {\sc Giovanni Scilla}\\
 Dipartimento di Matematica `G. Castelnuovo'\\
 `Sapienza' Universit\`a di Roma\\
 piazzale Aldo Moro 5, 00185 Roma - ITALY\\
 (scilla@mat.uniroma1.it)
                    \end{center}
                    \vspace{1cm}

\hspace*{-0.6cm}{\bf Abstract.} In this paper we describe the asymptotic behavior of rigid spin lattice energies by exhibiting a continuous interfacial limit energy as scaling to zero the lattice spacing. The limit is not trivial below a percolation threshold: it can be characterized by two phases separated by an interface. The macroscopic surface tension at this interface is defined through a first-passage percolation formula, related to the chemical distance on the lattice $\Z^2$. We also show a continuity result, that is the homogenization of rigid spin system is a limit case of the elliptic random homogenization.

\section{Introduction}
In the context of Variational theories in Materials Science it is often necessary to model media with fine microstructure and to describe their properties via averaged effective energies. This is the main goal of Homogenization theory (see e.g. \cite{BraDeF, ChePia}). In some cases periodic microstructure is not sufficient, so that random media have to be considered.

The model problem that we have in mind is that of a crystalline two\hbox{-}dimensional solid subject to fracture.
We suppose that the relevant scale is that of the surface (fracture) energy, so we may neglect the elastic energy of the lattice (this can be taken separately into account as in the paper \cite{BraPia}). In this case, depending on the applied forces or boundary displacement of the sample, a fracture may appear, separating two regions where the displacement is constant. In the Griffith theory of Fracture (see \cite{Gf}), the energy necessary for the creation of a crack is proportional to its area; in a discrete setting this is translated in the number of atomic bonds that are broken. In our model, at the atomistic level, there is a random distribution of `strong' unbreakable bonds and `weak' (ferromagnetic) breakable bonds. This model translates into a \emph{rigid} spin problem, where the two values of the spin parametrize the two regions of constant displacement of the crystal. We note that in this problem the random distribution of rigid or weak bonds is considered as fixed and as characteristic of the crystalline material, so that we are interested in almost sure properties of the overall energies when the measure of the sample is large with respect to the atomic distance.

The way we will describe the overall behavior of this system is by scaling the domain lattice by a small parameter $\epsilon$ and introducing the corresponding scaled energies, and then compute the variational limit ($\Gamma$\hbox{-}limit) of such energies, which is defined on the continuum and it can be considered as an effective energy.

The microscopic energy under examination can be written as

\begin{equation}
\sum_{ij}\sigma_{ij}^\omega (1-u_iu_j),
\label{system1}
\end{equation}
\\
where $u_i\in\{\pm1\}$ is a spin variable indexed on the lattice $\Z^2$, the sum runs on nearest neighbors (i.e. $|i-j|=1$) in a given portion $\Omega\cap\Z^2$ of $\Z^2$, the coefficients $\sigma_{ij}^\omega$ depend on the realization $\omega$ of an {independent and identically distributed (i.i.d.)} random variable and

\begin{equation}
\sigma_{ij}^\omega=
\begin{cases}
+\infty& \mbox{with probability}\quad p\\
1&\mbox{with probability}\quad 1-p,
\end{cases}
\end{equation}
\\
with $p\in[0,1]$ fixed and the convention $+\infty\cdot0=0$. 
In place of (\ref{system1}) we could consider the energies

\begin{equation*}
-\sum_{ij}\sigma_{ij}^\omega u_iu_j,
\label{system2}
\end{equation*}
\\
but in this case, just to avoid ambiguities in the sum, $\sigma_{ij}^\omega=+\infty$ forces $u_i=u_j$ and this gives a constraint for the problem.

In recent papers {\scshape Braides} and {\scshape Piatnitski} \cite{BraPia2,BraPia3} treated the cases of \emph{elliptic} random spin energies, that is with equi\hbox{-}bounded strictly positive random coefficients, and of \emph{dilute} spin energies, with random coefficients given by

\begin{equation*}
\widetilde{\sigma}_{ij}^\omega=
\begin{cases}
1& \mbox{with probability}\quad p\\
0&\mbox{with probability}\quad 1-p.
\end{cases}
\end{equation*}

In order to describe the behavior as the size of $\Omega$ diverges we introduce a scaled problem, as is customary in the passage from lattice systems to continuous variational problems, in which, on the contrary, ${\Omega}$ is kept fixed, but scaled energies are defined as follows. A small parameter $\epsilon>0$ is introduced, the lattice is scaled accordingly to $\epsilon\Z^2$, and the energies (\ref{system1}) are scaled (after multiplying by 2) to

\begin{equation}
E_\varepsilon^\omega(u):=\sum_{ij}\epsilon\sigma_{ij}^\omega (u_i-u_j)^2.
\label{rescaled}
\end{equation}
Note that uniform states (which are pointwise minimizers of the ``integrand'') have zero energy; moreover, the ``surface scaling'' $\epsilon$ is driven by the knowledge that for $p = 0$ (i.e., for ferromagnetic interactions) the $\Gamma$\hbox{-}limit with that scaling is not trivial (as shown e.g. by {\scshape Alicandro, Braides} and {\scshape Cicalese}~\cite{ABC}). After this scaling, the sum is taken on nearest neighbors in $\Omega\cap\epsilon\Z^2$, and the normalization allows also to consider $\Omega = \R^2$ (in this case the domain of the energy is composed of all $u$ which are constant outside a bounded set).

The coarse graining of these energies corresponds to a general approach in the theory of $\Gamma$\hbox{-}convergence for lattice systems where the discrete functions $u =\{u_i\}$ are identified with their piecewise-constant extensions, and the scaled lattice energies with energies on the continuum whose asymptotic behavior is described by taking $L^1$\hbox{-}limits in the $u$ variable and applying a mesoscopic homogenization process to the energies. A general theory for interfacial energies by {\scshape Ambrosio} and {\scshape Braides} \cite{AmbBra} suggests the identification of limit energies with functionals of the form

\begin{equation*}
\int_{\Omega\cap\partial\{u=1\}}\varphi(x,\nu)\,d\mathcal{H}^1,
\end{equation*}
with $\nu$ the normal to $\partial\{u=1\}$.

{Our analysis will be carried out by using results from Percolation theory. Percolation is a model for random media (see \cite{Grimmett, Kesten1}). We are interested in \emph{bond percolation} on the square lattice $\Z^2$: we view $\Z^2$ as a graph with edges between neighboring vertices, and all edges are, independently of each other, chosen to be `strong' with probability $p$ and `weak' with probability $1-p$. A weak path is a sequence of consecutive weak edges, a weak cluster is a maximal connected component of the collection of weak edges. Percolation exhibits a phase transition: there exists a critical value of probability $p_c$, the \emph{percolation threshold}, such that if $p<p_c$ then with probability one there exists a unique infinite weak cluster, while if $p>p_c$ then all the weak clusters are finite almost surely. For bond percolation on $\Z^2$, the percolation threshold is given by $p_c=\frac{1}{2}$.}

Actually, the structure of the $\Gamma$\hbox{-}limit of the energies (\ref{rescaled}) depends on probability through the percolation threshold. Above the percolation threshold the $\Gamma$\hbox{-}limit is $+\infty$ on the functions not identically equal to 1 or -1. Below the percolation threshold, instead, the coarse graining leads first to showing that indeed we may define a limit magnetization $u$ taking values in $\{\pm1\}$. This $u$ is obtained as a $L^1$\hbox{-}limit on the scaled infinite weak cluster, thus neglecting the values $u_i$ on nodes $i$ isolated from that cluster. The surface tension is obtained by optimizing the almost sure contribution of the interfaces, and showing that it can be expressed as a first\hbox{-}passage percolation problem, so that the limit is of the form

\begin{equation}
\int_{{\Omega}\cap\partial\{u=1\}}\lambda_p(\nu)\,d\mathcal{H}^1.
\label{limiten}
\end{equation}
\\
{The $\Gamma$\hbox{-}$\mathop{\lim\inf}$ inequality is obtained by a blow\hbox{-}up argument. We perform a construction based on the Channel property (Theorem~\ref{channel}) which allows to modify the test sets in order to get a `weak' boundary, avoiding bonds with infinite energy. This  is useful also for the construction of a recovery sequence.}

This type of variational percolation results can be linked to the paper by {\scshape Braides} and {\scshape Piatnitski} \cite{BraPia} where discrete fracture of a membrane is studied and linked to large deviations for the chemical distance in supercritical Bernoulli percolation. The value $\lambda_p(\nu)$ is defined through the asymptotic behavior of the chemical distance (that is, the distance on the infinite weak cluster) between a pair of points aligned with $\nu$. The general framework for first\hbox{-}passage percolation and chemical distance can be found in \cite{GK, Kesten3}. The result of this paper is that in the subcritical case, a crack in the crystal may appear following a minimal path on the infinite weak cluster and the microscopical pattern of the lattice (this fact justifies the anisotropy of the fracture energy (\ref{limiten})). In the supercritical case, instead, the solid almost surely is rigid and there is no fracture. 

The paper is organized as follows. {In Section~\ref{gmt} we fix notation and recall some definitions and results from Geometric Measure Theory and $\Gamma$\hbox{-}convergence in dimension $d\geq2$, even though the main result of this paper holds in $\R^2$.} In Section~\ref{setting} we present the setting of the problem describing the energies that we will consider. Section~\ref{recall} contains some results from percolation theory necessary for the computations. Section~\ref{percolationthm} contains the proof of the main result: in the subcritical regime (that is, {for  $p<p_c$}), the energies $\Gamma$\hbox{-}converge to a deterministic anisotropic perimeter whose density is obtained by means of an asymptotic first\hbox{-}passage percolation formula related to the chemical distance on the lattice; in the supercritical regime (that is, {for $p>p_c$}), the $\Gamma$\hbox{-}limit is identically $+\infty$. In Section~\ref{limitcase} we show that the homogenization of rigid spin systems is actually a limit case of the elliptic random homogenization of spin systems; that is, the behavior of a rigid spin system is approximated by that of an elliptic spin system with one of the interaction coefficients very large. The proof of this new ``continuity'' property of the surface tension (Proposition~\ref{continuo}) essentially relies on a percolation result (Lemma~\ref{percentage}).

\section{Notation}\label{gmt}
Let $d\geq2$. If $A$ is a measurable subset of $\R^d$, we denote its $d$\hbox{-}dimensional Lebesgue measure indifferently by $\mathcal{L}^d(A)$ or $|A|$. $\mathcal{H}^k$ denotes the $k$\hbox{-}dimensional Hausdorff measure. $B_\rho(x)$ is the open ball of center $x$ and radius $\rho$ and $S^{d-1}$ is the boundary of $B_1(0)$. If $\nu\in S^1$, $Q_\rho^\nu(x)$ is the square centered at $x$, of side length $\rho$ and one side orthogonal to $\nu$, that is

\begin{equation*}
Q_\rho^\nu(x)=\{y\in\R^2:|\langle y-x,\nu\rangle|\leq \rho/2, |\langle y-x,\nu^\perp\rangle|\leq \rho/2\},
\end{equation*} 
\\
where $\nu^\perp=(-\nu_2,\nu_1)$ denotes the clockwise rotation by $\pi/2$ of $\nu$.

\subsection{Functions of bounded variation and sets of finite perimeter}\label{richiami}

For the general theory of functions of bounded variation and sets of finite perimeter we refer to \cite{AFP,Bra98}; we recall some definitions and results necessary in the sequel. Let $\Omega$ be an open subset of $\R^d$. We say that $u\in L^1({\Omega})$ is a \emph{function of bounded variation} if its distributional first derivatives $D_i u$ are Radon measures with finite total variation in $\Omega$. We denote this space by $BV({\Omega})$ and we write $u\in BV({\Omega};\{\pm1\})$ when the function $u$ is of bounded variation in ${\Omega}$ and takes only the values -1 and +1.

Let $u:{\Omega}\longrightarrow\R$ be a Borel function. We say that $z\in\R$ is the \emph{approximate limit} of $u$ {at} $x$ if for every $\epsilon>0$

\begin{equation*}
\lim_{\rho\to0^+}\rho^{-d}\mathcal{L}^d(\{y\in B_\rho(x)\cap{\Omega}:|u(y)-z|>\epsilon\})=0.
\end{equation*}
\\
We define the \emph{jump set $S(u)$ of function $u$} as the subset of ${\Omega}$ where the approximate limit of $u$ does not exist. It turns out that $S(u)$ is a Borel set and $\mathcal{L}^d(S(u))=0$. If $u\in BV({\Omega})$, then $S(u)$ is \emph{countably $(d-1)$\hbox{-}rectifiable}; that is, $S(u)=N\cup\left(\bigcup_{i\in\N}K_i\right)$, where $\mathcal{H}^{d-1}(N)=0$ and $(K_i)$ is a sequence of compact sets, each contained in a $C^1$ hypersurface $\Gamma_i$. A \emph{normal unit vector} $\nu_u$ to $S(u)$ exists $\mathcal{H}^{d-1}$\hbox{-}almost everywhere on $S(u)$, in the sense that, if $S(u)$ is represented as above, then $\nu_u(x)$ is normal to $\Gamma_i$ for $\mathcal{H}^{d-1}$\hbox{-}almost everywhere $x\in K_i$.

Let $E$ be a Borel subset of $\R^d$. The \emph{essential boundary} $\partial^*E$ of $E$ is defined as

\begin{equation*}
\partial^*E=\left\{x\in\R^d:\mathop{\lim\sup}_{\rho\to0}\frac{\mathcal{L}^d(B_\rho(x)\cap E)}{\rho^d}>0, \mathop{\lim\sup}_{\rho\to0}\frac{\mathcal{L}^d(B_\rho(x)\backslash E)}{\rho^d}>0\right\}.
\end{equation*}
\\
The set $E$ is of \emph{finite perimeter} in ${\Omega}$ if its characteristic function $\chi_E$ is in $BV({\Omega})$. The total variation $|D\chi_E|({\Omega})$ is the \emph{perimeter} of $E$ in ${\Omega}$, denoted by $P(E;{\Omega})$. For $\mathcal{H}^{d-1}$\hbox{-}almost every $x\in\partial^*E$, the limit

\begin{equation*}
\nu_E(x)=\lim_{\rho\to0}\frac{D\chi_E(B_\rho(x))}{|D\chi_E|(B_\rho(x))}
\end{equation*}
\\
exists and belongs to $S^{d-1}$; the vector $\nu_E$ is the generalized \emph{inner normal} to $\partial^*E$. The set of points $x\in\text{supp}(D\chi_E)\cap{\Omega}$ where this property holds is called the \emph{reduced boundary} of $E$. For any $x$ in the reduced boundary of $E$, the sets $(E-x)/\rho$ locally converge in measure as $\rho\to0$ to the half\hbox{-}space orthogonal to $\nu_E(x)$ and containing $\nu_E(x)$. The measure $D\chi_E$ can be represented as

\begin{equation*}
D\chi_E=\nu_E\mathcal{H}^{d-1}\res\partial^*E.
\end{equation*}
\\
In particular, for every set $E$ of finite perimeter in ${\Omega}$, $P(E;{\Omega})=\mathcal{H}^{d-1}(\partial^*E\cap{\Omega})$.

\subsection{$\Gamma$\hbox{-}convergence}

For the definition and properties of $\Gamma$\hbox{-}convergence we refer to \cite{Bra02,Bra06,BraDeF,DalM}. We just recall the definition of \emph{$\Gamma$\hbox{-}convergence} of a family of functionals $(F_\epsilon)_{\epsilon>0}$ defined on $BV({\Omega})$: we say that $(F_\epsilon)$ $\Gamma$\hbox{-}converges to $F$ (on $BV({\Omega})$ with respect to the convergence in measure) if for all $u\in BV({\Omega})$ and for all sequences $(\epsilon_j)$ of positive numbers converging to 0

\begin{description}
\item[(i)] (\emph{lower bound}) for all sequences $(u_{\epsilon_j})$ converging to $u$ in measure we have

\begin{equation*}
F(u)\leq\mathop{\lim\inf}_{j\to+\infty}F_{\epsilon_j}(u_{\epsilon_j});
\end{equation*}
\item[(ii)] (\emph{upper bound}) there exists a sequence $(u_{\epsilon_j})$ converging to $u$ in measure such that

\begin{equation*}
F(u)\geq\mathop{\lim\sup}_{j\to+\infty}F_{\epsilon_j}(u_{\epsilon_j}).
\end{equation*}
\end{description} 
If {\bf(i)} and {\bf(ii)} hold then we write $F(u)=\Gamma\hbox{-}\displaystyle\lim_{\epsilon\to0}F_\epsilon(u)$. We define the \emph{$\Gamma$\hbox{-}lower limit} as

\begin{equation*}
\Gamma\hbox{-}\mathop{\lim\inf}_{\epsilon\to0}F_\epsilon(u)=\inf\left\{\mathop{\lim\inf}_{j\to+\infty}F_{\epsilon_j}(u_{\epsilon_j}):u_{\epsilon_j}\to u\right\}
\end{equation*}
\\
and the \emph{$\Gamma$\hbox{-}upper limit} as

\begin{equation*}
\Gamma\hbox{-}\mathop{\lim\sup}_{\epsilon\to0}F_\epsilon(u)=\inf\left\{\mathop{\lim\sup}_{j\to+\infty}F_{\epsilon_j}(u_{\epsilon_j}):u_{\epsilon_j}\to u\right\},
\end{equation*}
\\
respectively. 

Then {\bf(i)} also reads $F(u)\leq\Gamma\hbox{-}\displaystyle\mathop{\lim\inf}_{\epsilon\to0}F_\epsilon(u)$ and {\bf(ii)} reads $F(u)\geq\Gamma\hbox{-}\displaystyle\mathop{\lim\sup}_{\epsilon\to0}F_\epsilon(u)$.

\section{Setting of the problem}\label{setting}
We use the notation for \emph{bond\hbox{-}percolation} problems as in {\scshape Braides-Piatnitski}~\cite{BraPia}, Section~2.4. In this percolation model, {we assign the label ``strong'' or ``weak'' to a bond with probability $p$ and $1-p$, respectively,} the choice being independent on distinct bonds.

Denote by $\hat{\Z}^2$ the \emph{dual grid} of $\Z^2$, that is the set of the middle points of the segments $[i,j]$, $i,j\in\Z^2$, $|i-j|=1$, of the standard integer grid $\Z^2$. The notation $i(\hat{z}),j(\hat{z})$ is used for the endpoints of the segment containing $\hat{z}$. We may identify each point $\hat{z}\in\hat{\Z}^2$ with the corresponding closed segment $[i(\hat{z}),j(\hat{z})]$, so that points in $\hat{\Z}^2$ are identified with bonds in $\Z^2$.

Let $({\Sigma},\mathcal{F},\textbf{P})$ be a probability space, and let $\{\xi_{\hat{z}},\hat{z}\in\hat{\Z}^2\}$ be a family of independent identically distributed random variables such that

$$
\xi_{\hat{z}}=
\begin{cases}
1& \mbox{(``strong'') with probability}\quad p\\
0&\mbox{(``weak'') with probability}\quad 1-p,
\end{cases}
$$
and $p\in[0,1]$ is fixed.
Let $\omega\in\Sigma$ be a realization of this i.i.d. random variable in $\Z^2$ and introduce the coefficients

$$
\sigma^\omega_{\hat{z}}=
\begin{cases}
+\infty& \mbox{if}\quad\xi_{\hat{z}}(\omega)=1\\
1&\mbox{otherwise}.
\end{cases}
$$
We write $\sigma^\omega_{\hat{z}}=\sigma^\omega_{ij}$, after identifying each $\hat{z}$ with a pair of nearest neighbors in $\Z^2$.

For each $\omega$ we consider the energies
\begin{equation}
E_\epsilon^\omega(u)=\frac{1}{8}\sum_{i,j\in \Omega_\epsilon}\epsilon\sigma^\omega_{ij}(u_i-u_j)^2,\quad\mbox{(with the convention $\infty\cdot0=0$)}
\label{energy}
\end{equation}
defined on $u:\Omega_\epsilon\to\{\pm 1\}$, where we use the notation $\Omega_\epsilon=\frac{1}{\epsilon}\Omega\cap \Z^2$ and $\Omega$ is an open subset of $\R^2$ with Lipschitz boundary. The factor $8$ is a normalization factor due to the fact that each bond is accounted {twice} and $(u_i-u_j)^2\in\{0,4\}$.

Each function $u:\Omega_\epsilon\to\{\pm 1\}$ is identified with the piecewise\hbox{-}constant function {(which, with a slight abuse of notation, we also denote by $u$)} such that $u(x)=u_i$ on each coordinate square of center $\epsilon i$ and side length $\epsilon$ contained in $\Omega$ and 1 otherwise, {no matter what this value is.} In this way all $u$ can be considered as functions in $L^1(\Omega)$, and more precisely in $BV(\Omega;\{\pm 1\})$.

The case $p=0$ corresponds to a \emph{ferromagnetic spin system}, which can be described approximately as $\epsilon\to0$ by the anisotropic perimeter energy (see \cite{ABC})

\begin{equation*}
F_0(u)=\int_{\partial^*\{u=1\}\cap{\Omega}}\|\nu_u\|_1\,d\mathcal{H}^1
\end{equation*}
defined on $u\in BV(\Omega;\{\pm 1\})$; here $\partial^*\{u=1\}$ denotes the measure\hbox{-}theoretical reduced boundary of the set of finite perimeter $\{u=1\}$, $\nu_u$ its inner normal (see Subsection~\ref{richiami}) and $\|\cdot\|_1$ the $\ell^1$\hbox{-}norm defined by

\begin{equation*}
\|x\|_1=|x_1|+|x_2|,
\end{equation*}
\\
{where $x=(x_1,x_2)$.}

{In this passage from discrete to continuous}, we identify each function $u:\Omega_\epsilon\to\{\pm 1\}$ with the set $A=\bigcup\{\epsilon i+\epsilon Q:i\in{\Omega}_\epsilon:u_i=1\}$ or the function $u\in BV(\Omega;\{\pm 1\})$ given by $u=2\chi_A-1$, where {$Q=\left[-\frac{1}{2},\frac{1}{2}\right)\times\left[-\frac{1}{2},\frac{1}{2}\right)$ denotes the coordinate semi\hbox{-}open unit square in $\R^2$ centered at 0.}

\section{Some results from percolation theory}\label{recall}
We recall some results from percolation theory (see \cite{BraPia}, and {\scshape Grimmett}~\cite{Grimmett}, {\scshape Kesten}~\cite{Kesten1} for general references on percolation theory).

We introduce a terminology for {\emph{weak points}}, that is those points $\hat{z}\in\hat{\Z}^2$ such that $\xi_{\hat{z}}=~0$. {Keeping in mind the identification of $\hat{z}$ with $[i(\hat{z}),j(\hat{z})]$ stated in the previous section, we denote the corresponding bond also by $\hat{z}$, and we refer to $\hat{z}$ as a \emph{weak bond} or point, indifferently.} We say that two weak points $\hat{z}$ and $\hat{z}'$ are \emph{adjacent} if the corresponding two segments have an endpoint in common. A sequence of weak bonds $\gamma=\{\hat{z}_1,\dots,\hat{z}_k\}$ is said to be a \emph{weak path} if any two consecutive points of this sequence are adjacent. In what follows we identify a weak path with the subset of $\R^2$ composed of the union of the corresponding segments; the \emph{length} of the weak path $\gamma$ is the number of its connections, and we denote it by $|\gamma|$. A subset $A$ of $\hat{\Z}^2$ of weak points is said to be \emph{connected} if for every two points $\hat{z}',\hat{z}''$ of $A$ there exists a weak path as above such that {$\forall j\in\{1,2,\dots,k\}$, $\hat{z}_j\in A,\hat{z}_1=\hat{z}',\hat{z}_k=\hat{z}''$.} A maximum connected component of adjacent weak points is called a \emph{weak cluster}.

A first result deals with the existence of infinite weak clusters and it shows that there is a critical probability ($p_c=\frac{1}{2}$ in two dimensions for the square lattice) which separates two different behaviors of the bond\hbox{-}percolation system. 

\begin{thm}[Percolation threshold]
For any $p>p_c=1/2$ \emph{(supercritical regime)} all the weak clusters are almost surely finite, while for any $p<1/2$ \emph{(subcritical regime)} with probability one there is exactly one infinite weak cluster $\mathcal{W}^\omega$.
\end{thm}
From now on we will refer to $\mathcal{W}^\omega$ simply as \emph{the} weak cluster.
Let $\nu=(\nu_1,\nu_2)\in S^1$ and $0<\delta<1$. We denote by $T_\nu^\delta$ the rectangle
\begin{equation*}
T_\nu^\delta=\{x\in\R^2:|\langle x,\nu\rangle|\leq\delta/2, 0\leq\langle x,\nu^\perp\rangle\leq1\}.
\end{equation*}
A path joining the smaller sides of the rectangle will be called a \emph{channel} (or \emph{left\hbox{-}right crossing}). A weak path with this property will be called a \emph{weak channel}. The following result gives a lower bound on the number of weak channels almost surely crossing sufficiently large rectangles (or squares) in the subcritical regime.

\begin{thm}[Channel property]
\label{channel}
Let $p<1/2$, {$\omega$ denote a realization} and $M>0$. Then there exist constants $c(p)>0$ and $c_1(p)>0$ such that almost surely for any $\delta$, $0<\delta\leq1$ {there is a number} $N_0=N_0(\omega,\delta)$ such that for all $N>N_0$ and for any $T_\nu^\delta$ and $|x_0|\leq M$ the rectangle $N(T_\nu^\delta+x_0)$ contains at least $c(p)N\delta$ disjoint weak channels which connect the smaller sides of the rectangle. Moreover, the length of each such a channel does not exceed $c_1(p)N$.
\end{thm}

{A realization $\omega\in\Sigma$ is said to be \emph{typical} if the statement of the Theorem~\ref{channel} holds for such an $\omega$.} Now we introduce some terminology also for \emph{strong bonds}, that is those points $\hat{z}\in\hat{\Z}^2$ such that $\xi_{\hat{z}}=1$. We consider the shifted lattice $Z_b=\Z^2+\left(\frac{1}{2},\frac{1}{2}\right)$ and notice that the set of middle points of its bonds coincides with $\hat{\Z}^2$. Thus, to each point $\hat{z}\in\hat{\Z}^2$ we can associate the corresponding bond in $Z_b$. If $\hat{z}$ is identified with the corresponding segment with endpoints in $Z_b$, then we may define the notion of \emph{adjacent points} as for weak bonds. The notion of a \emph{strong channel} and a \emph{strong cluster} is modified accordingly. For $p>1/2$ there is almost surely a unique infinite strong cluster and the channel property stated above holds for the strong channels as well.

To simplify the notation, from now on we will denote by $x$ (in place of $\hat{x}$) a generic point in $\hat{\Z}^2$. If $p<1/2$ and two points $x,y\in\hat{\Z}^2$ belong to the weak cluster, then by definition of cluster there is at least a path $\gamma$ in the cluster joining $x$ and $y$. To describe the metric properties of the weak cluster we introduce a random distance.
\begin{defn}
Let $x,y\in{\hat{\Z}}^2$ and $\omega$ be a realization of the random variable. The \emph{chemical distance} $D^\omega(x,y)$ between $x$ and $y$ in the realization $\omega$ is defined as
\begin{equation}
D^\omega(x,y)=\min_\gamma|\gamma|,
\end{equation}
where $|\gamma|$ is the length of the path $\gamma$ and the minimum is taken on the set of paths joining $x$ and $y$ and that are weak in the realization $\omega$.
\end{defn}

\begin{oss}
The chemical distance is defined only if $x$ and $y$ are in the same cluster; otherwise, by convention, $D^\omega(x,y)=+\infty$.
\end{oss}

{When it is finite, the random distance} $D^\omega(x,y)$ is thus the minimal number of weak bonds needed to link $x$ and $y$ in the realization $\omega$ (also $x$ and $y$ are taken {into} account), and is thus {not smaller} than $\|x-y\|_1$.
When $p<\frac{1}{2}$, $D^\omega(0,x)$ on {the weak cluster} can be seen as a \emph{travel time} between $0$ and $x$ in a first\hbox{-}passage percolation model (see~\cite{GM,Kesten2}) where the passage times of the edges are independent identically distributed random variables with common distribution $p\delta_{+\infty}+(1-p)\delta_1$. The following Lemma deals with the existence of an asymptotic time constant in a given direction.
\begin{lem}
Assume that $0\in\mathcal{W}^\omega$. For any $\tau\in\R^2$ the following limit exists almost surely and does not depend on $\omega$
\begin{equation}
\lambda_p(\tau)=\lim_{\substack{m\to +\infty \\ 
                0\longleftrightarrow\lfloor m\tau\rfloor }}
\frac{1}{m}D^\omega(0,\lfloor m\tau\rfloor),
\label{limit}
\end{equation}
where {$\lfloor m\tau\rfloor=(\lfloor m\tau\rfloor_1,\lfloor m\tau\rfloor_2)$}, $\lfloor m\tau\rfloor_k=\lfloor m\tau_k\rfloor$ is the integer part of the $k$-th component of $m\tau$ and $0\longleftrightarrow\lfloor m\tau\rfloor$ means that $0$ and $\lfloor m\tau\rfloor$ are linked by a path in the weak cluster. Moreover, $\lambda_p$ defines a norm on ${\mathbb R}^2$.
\label{lim}
\end{lem}
\proof
See {\scshape Garet-Marchand~\cite{GM1}}, Corollary 3.3.
\endproof

\bigskip
{The same asymptotic result holds for sequences of points in the weak cluster `asymptotically aligned' with $\tau$, that is $x_m,y_m$ such that $y_m-x_m=m\tau+o(m)$ as $m\to\infty$. The proof of this fact relies essentially on the following large deviation result for the chemical distance (see {\scshape Garet-Marchand}~\cite{GM}):}

\begin{thm}\label{deviation}
{Let $p<1/2$ and $\lambda_p$ be the norm on $\R^2$ given by Lemma~\ref{lim}. Then
\begin{equation}
\forall \epsilon>0,\quad \mathop{\lim\sup}_{m\to+\infty}\frac{\log {\bf  P}[0\longleftrightarrow\lfloor m\tau\rfloor ,D^\omega(0,\lfloor m\tau\rfloor)/\lambda_p(\tau)\notin(1-\epsilon,1+\epsilon)]}{m}<0.
\label{ld}
\end{equation}
}
\end{thm} 

\begin{prop}
Let $(x_m), (y_m)$ be two sequences of points in $\hat{\Z}^2$ contained in the weak cluster such that 
\begin{equation}
\sup_m\left\{\frac{|x_m|}{m}+\frac{|y_m|}{m}\right\}\leq C<+\infty\quad\text{and}\quad y_m-x_m=m\tau+o(m),
\label{asymp}
\end{equation}
where $\tau\in\R^2$, $C$ is a positive constant and $o(m)/m\to0$ as $m\to +\infty$. Then the following limit exists almost surely and does not depend on $\omega$
\begin{equation}
\lambda_p(\tau)=\lim_{m\to +\infty}
\frac{1}{m}D^\omega(x_m,y_m).
\label{samelim}
\end{equation}
\label{generallim}
\end{prop}
\proof 
We denote by $\widetilde{\lambda}_p(\tau)$ the right hand side of equation (\ref{samelim}). We prove (\ref{samelim}) in the case that $y_m=\lfloor mx+m\tau\rfloor$ and $x_m=~\lfloor mx\rfloor$, $x\neq0$. The stationarity of the i.i.d. Bernoulli process ensures that the probability law of $D^\omega(\lfloor mx\rfloor, \lfloor mx+~m\tau\rfloor)$ is the same of $D^\omega(0, \lfloor m\tau\rfloor)$.

Therefore, by (\ref{ld}) we have that

\begin{equation*}
\forall \epsilon>0,\quad \mathop{\lim\sup}_{m\to+\infty}\frac{\log {\bf  P}\left[\lfloor mx\rfloor\longleftrightarrow\lfloor mx+m\tau\rfloor ,\frac{D^\omega(\lfloor mx\rfloor,\lfloor mx+m\tau\rfloor)}{\lambda_p(\tau)}\notin(1-\epsilon,1+\epsilon)\right]}{m}<0.
\label{ld2}
\end{equation*}
By Borel\hbox{-}Cantelli Lemma we obtain that $\forall \epsilon>0$, 

\begin{equation*}
\mathop{\lim\sup}_{\substack{m\to +\infty \\ 
           \lfloor mx\rfloor\longleftrightarrow\lfloor mx+m\tau\rfloor }}
\frac{1}{m}D^\omega(\lfloor mx\rfloor,\lfloor mx+m\tau\rfloor)\in[(1-\epsilon)\lambda_p(\tau),(1+\epsilon)\lambda_p(\tau)]
\end{equation*}
{\bf P}\hbox{-}almost surely, that is, 

\begin{equation*}
\lim_{\substack{m\to +\infty \\ 
           \lfloor mx\rfloor\longleftrightarrow\lfloor mx+m\tau\rfloor }}
\frac{1}{m}D^\omega(\lfloor mx\rfloor,\lfloor mx+m\tau\rfloor)=\lambda_p(\tau)
\end{equation*}
{\bf P}\hbox{-}almost surely.
{By a compactness argument, we have that $\widetilde{\lambda}_p(\tau)\leq\lambda_p(\tau)$. Indeed, if $x_m,y_m$ are as in (\ref{asymp}), then there exist a subsequence $m_j\to+\infty$ and $\widetilde{x},\widetilde{y}$ such that $\displaystyle\frac{x_{m_j}}{m_j}\to\widetilde{x},\displaystyle\frac{y_{m_j}}{m_j}\to\widetilde{y}$, with $\widetilde{y}=\widetilde{x}+\tau$. Therefore, in the computation of the limit in (\ref{samelim}) we can choose $x_m=\lfloor m\widetilde{x}\rfloor$ and $y_m=\lfloor m\widetilde{x}+m\tau\rfloor$.} 

Now, if we consider two points $x_m,y_m$ on the weak cluster satisfying (\ref{asymp}), we can find $x$ such that $x_m$ and $\lfloor mx\rfloor$, $y_m$ and $\lfloor mx+m\tau\rfloor$ almost surely are linked by weak paths whose length is at most $o(m)$. Hence,

\begin{equation*}
D^\omega(\lfloor mx\rfloor,\lfloor mx+m\tau\rfloor)\leq D^\omega(x_m,y_m)+o(m)
\end{equation*}
\\
and from this it follows that $\lambda_p(\tau)\leq\widetilde{\lambda}_p(\tau)$.
\endproof

\begin{oss}
If $\nu\in\R^2$ is a unit vector and $\tau=\nu^\perp$, then by symmetry $\lambda_p(\nu)=\lambda_p(\tau)$.
\label{symmetry}
\end{oss}

\section{The rigid percolation theorem}\label{percolationthm}

We first remark that the energies $E^\omega_\epsilon$ defined by (\ref{energy}) are equi\hbox{-}coercive with respect to the strong $L^1$\hbox{-}convergence. {The proof is immediate as in the elliptic case, while for dilute spin energies this result requires a more difficult argument (see \cite{BraPia3}, Section 3.1).}

\begin{lem}[Equi\hbox{-}coerciveness of $E_\epsilon^\omega$]
Let $\Omega$ be a bounded Lipschitz open set. {For any $\omega$ in a set $\widetilde{\Sigma}\subseteq\Sigma$ with ${\mathbf P}(\widetilde{\Sigma})=1$,} if $(u_{\epsilon_k})$ is a sequence such that $\sup_kE_{\epsilon_k}^\omega(u_{\epsilon_k})<+\infty$, then there exists a function $u\in BV(\Omega;\{\pm1\})$ and a subsequence, still denoted by $(u_{\epsilon_k})$, such that $u_{\epsilon_k}\to u$ in $L^1({\Omega})$.
\end{lem}
\proof
Equi\hbox{-}boundedness of the energies forces the coefficients $\sigma_{ij}^\omega$ to be equal to $1$ almost surely if $(u_{\epsilon_k})_i\not=(u_{\epsilon_k})_j$, so that the equi\hbox{-}coerciveness follows from that of ferromagnetic energies (see e.g. \cite{BraPia2}, Section~2).
\endproof

\bigskip
The main result of this paper is the following.

\begin{thm}[Rigid percolation theorem]
Let $\Omega$ be a bounded Lipschitz open set and $E_\epsilon^\omega$ be the energies defined by (\ref{energy}). Then we have two regimes:
\begin{description}
\item[(a)] If $p<1/2$ \emph{(subcritical regime)}, then ${\mathbf P}$\hbox{-}almost surely there exists the $\Gamma$\hbox{-}limit of $E_\epsilon^\omega$ with respect to the $L^1({\Omega})$\hbox{-}convergence on $BV(\Omega;\{\pm1\})$, it is deterministic, and is given by
\begin{equation}
\label{energ}
F_p(u)=\int_{\Omega\cap\partial^*\{u=1\}}\lambda_p(\nu)\,d\mathcal{H}^1,
\end{equation}
for $u\in BV(\Omega;\{\pm1\})$. {In (\ref{energ}) $\lambda_p$ is defined by (\ref{limit}), (\ref{samelim})} and $\nu$ is the unit normal vector at $\partial^*\{u=1\}$.
\item[(b)] If $p>1/2$ \emph{(supercritical regime)}, then {\bf P}\hbox{-}almost surely there exists the $\Gamma$\hbox{-}limit of $E_\epsilon^\omega$ and it coincides with $F(u)\equiv+\infty$ on the whole $L^1(\Omega)$ except for $u$ constant identically $\pm 1$.
\end{description}
\label{percthm}
\end{thm}

\proof
{\bf (a)} We begin with the proof of the lower bound (liminf inequality), and fix a typical realization $\omega$ and a family $u_\epsilon\to u$  in $L^1({\Omega})$ with $u\in BV(\Omega;\{\pm1\})$ such that $\displaystyle\mathop{\lim\inf}_{\epsilon\to 0}E^\omega_\epsilon(u_\epsilon)<\infty$. We can assume, up to a subsequence, that such a liminf is actually a limit.

For all $\epsilon$ we consider the set in the dual lattice $\epsilon{\hat\Z}^2$ of $\epsilon\mathbb Z^2$ defined by

\begin{equation*}
S_\epsilon =\left\{\epsilon k: k=\frac{i+j}{2}, {i, j\in \Omega_\epsilon},|i-j|=1, u_\epsilon(\epsilon i)=1, u_\epsilon(\epsilon j)=-1\right\}
\end{equation*}
and the measure
\begin{equation*}
\mu_\epsilon=\sum_{\epsilon k\in S_\epsilon}\epsilon\sigma^\omega_k\delta_{\epsilon k}.
\end{equation*} 
Note that $E_\epsilon^\omega(u_\epsilon)=\mu_\epsilon(\Omega)$, so that the family of measures $\{\mu_\epsilon\}$ is equi-bounded. Hence, up to further subsequences, we can assume that $\mu_\epsilon\rightharpoonup^*\mu$, where $\mu$ is a finite measure.

Let $A=\{u=1\}$ and $A_\epsilon=\{u_\epsilon=1\}$. With fixed $h\in\mathbb N$, we consider the collection $\mathcal Q_h$ of squares $Q_\rho^\nu(x)$ 
such that the following conditions are satisfied:
\begin{description}
\item[(i)] $x\in\partial^*A\text{ and }\nu=\nu(x)$;
\item[(ii)]$\left|(Q_\rho^\nu(x)\cap A)\triangle \Pi^\nu(x)\right|\leq\frac{1}{h}\rho^2$, where $\Pi^\nu(x)=\{y\in\R^2:\langle y-x,\nu\rangle\geq0\}$;
\item[(iii)]$\displaystyle\left|\frac{\mu(Q_\rho^\nu(x))}{\rho}-\frac{d\mu}{d\mathcal{H}^1\res\partial^*A}(x)\right|\leq\frac{1}{h}$;
\item[(iv)]$\displaystyle\left|\frac{1}{\rho}\int_{Q_\rho^\nu(x)\cap\partial^*A}\lambda_p(\nu(y))\,d\mathcal{H}^1(y)-\lambda_p(\nu(x))\right|\leq\frac{1}{h}$;
\item[(v)] $\mu(Q_\rho^\nu(x))=\mu(\overline{Q_\rho^\nu(x)})$.
\end{description}
For fixed $x\in\partial^*A$ and for $\rho$ small enough {\bf (ii)} is satisfied by definition of reduced boundary (see subsection \ref{richiami}), {\bf(iii)} follows from the Besicovitch Derivation Theorem provided that
\begin{equation*}
\frac{d\mu}{d\mathcal{H}^1\res\partial^*A}(x)<+\infty;
\end{equation*}
{\bf(iv)} holds by the same reason ($x$ is a Lebesgue point of $\lambda_p$), and {\bf(v)} is satisfied for almost all $\rho>0$ since $\mu$ is a finite measure (and so $\mu(\partial Q^\nu_\rho(x))=0$).

We deduce that $\mathcal{Q}_h$ is a fine covering of the set
\begin{equation*}
\partial^*A_\mu=\left\{x\in\partial^*A:\frac{d\mu}{d\mathcal{H}^1\res\partial^*A}(x)<+\infty\right\},
\end{equation*}
so that (by Morse's lemma, see \cite{Morse}) there exists a countable family of disjoint closed squares $\{\overline{Q_{\rho_j}^{\nu_j}(x_j)}\}$ still covering $\partial^*A_\mu$. {Note that 
\begin{equation*}
\mathcal{H}^1(\partial^*A\backslash\partial^*A_\mu)=0
\end{equation*}
because of the existence of the derivative of the measure $\mu$ with respect to $\mathcal{H}^1\res\partial^*A$}.

We now fix one of such squares $Q_\rho^\nu(x)$. We would like to use the sets $\frac{1}{\epsilon}A_\epsilon$ as test sets to estimate from below the part of $\mu$ concentrated on $\partial^*A$. Since these sets could not have correct boundary data, we modify them on the boundary. {As in the case of spin energies with bounded coefficients (see \cite{BraPia2,BraPia3}),} we could truncate the sets with the hyperplane $\partial\Pi^\nu=\{\langle y-x,\nu\rangle=0\}$, but it may have infinite energy (possibly containing some strong bond), so we approximate it with a weak path $\gamma_\omega$. 

\begin{figure}
\centering
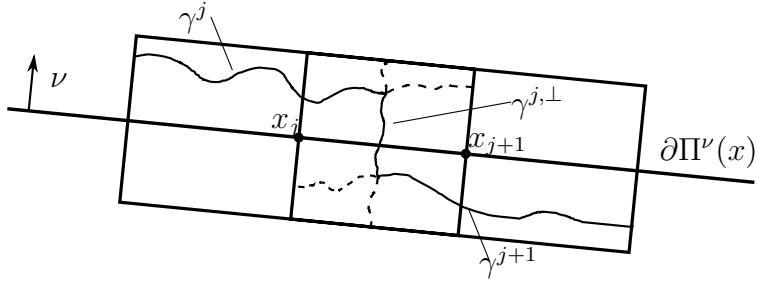
\caption{{Construction of $\widetilde{\gamma}^{j,j+1}$.}}
\label{construction}
\end{figure}

{We subdivide the construction of test sets into steps.}
\\
{{\bf Step 1: Construction of the weak path $\gamma_\omega$.}} Let $0<\delta\leq1$. {We cover the set $\partial\Pi^\nu(x)\cap Q_\rho^\nu(x)$ by considering the points} 
\begin{equation}
x_j=x+j\rho\delta\nu^\perp,\quad |j|=0,1,\dots,\left\lfloor\frac{1}{2\delta}\right\rfloor+1,
\label{points}
\end{equation}
{and the rectangles $R_j$ centered at $x_{j}$}, with side-lengths $\rho\delta$ and $2\rho\delta$ and the small sides parallel to $\nu$. The rectangles $R_j$ and $R_{j+1}$ have in common a square of side\hbox{-}length $\rho\delta$ and we denote it by $Q_{\rho\delta}^j$. From the channel property (Theorem~\ref{channel}) for $\epsilon$ small enough we can find a weak channel $\gamma^j$ (the highest) in $R_j$ and a weak channel $\gamma^{j+1}$ (the lowest) in $R_{j+1}$ whose length is at most $2c_1(p)\rho\delta/\epsilon$. Applying the same property to $Q_{\rho\delta}^j$ we can find a weak channel $\gamma^{j,\perp}$ connecting the two opposite sides of $Q_{\rho\delta}^j$ orthogonal to $\nu$, whose length is at most $c_2(p)\rho\delta/\epsilon$. The union $\gamma^j\cup\gamma^{j,\perp}\cup\gamma^{j+1}$ contains a weak path $\widetilde{\gamma}^{j,j+1}$ connecting the smaller sides of the rectangle $R_j\cup R_{j+1}$ (see Fig.~\ref{construction}). 

\noindent
{We can repeat this construction for $R_{j+1},R_{j+2}$ and $Q^{j+1}_{\rho\delta}$ choosing $\gamma^{j+1}$, the highest weak channel $\gamma^{j+2}$ in $R^{j+2}$ and $\gamma^{{j+1},\perp}$ in $Q^{j+1}_{\rho\delta}$ to define the weak path $\widetilde{\gamma}^{j+1,j+2}$. If we iterate this procedure} for any couple of rectangles $R_j,R_{j+1}$, $j$ as in (\ref{points}), the desired $\gamma_\omega$ will be obtained by gluing all the paths $\widetilde{\gamma}^{j,j+1}$.
\\
{{\bf Step 2: Estimates.}} Note that $\gamma_\omega$ disconnects $Q_\rho^\nu(x)$ and we denote by $Q^+_\omega$ the connected component of $Q_\rho^\nu(x)$ containing $Q_\rho^\nu(x)\cap\{\langle y-x,\nu\rangle\geq\rho\delta/2\}$.
We have that
\begin{equation}
\left|(Q_\rho^\nu(x)\cap \Pi^\nu(x))\triangle Q^+_\omega)\right|\leq\rho^2\delta.
\label{estimate}
\end{equation}

\begin{figure}
\centering
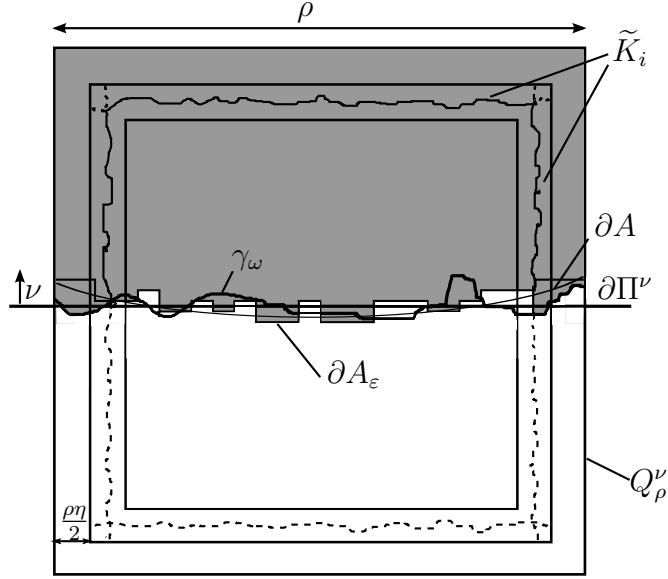
\caption{{Construction of a test set.}}
\label{testset}
\end{figure}

By {\bf(ii)} and (\ref{estimate}), choosing $\epsilon$ small enough and using the fact that $|A_\epsilon\triangle A|\to 0$, we obtain
\begin{equation}
\begin{split}\label{(16)}
\left|(Q_\rho^\nu(x)\cap A_\epsilon)\triangle Q^+_\omega\right|&\leq\left|(Q_\rho^\nu(x)\cap A)\triangle \Pi^\nu(x)\right|+\left|Q_\rho^\nu(x)\cap (A_\epsilon\triangle A)\right|+\\
&+\left|(Q_\rho^\nu(x)\cap \Pi^\nu(x))\triangle Q^+_\omega)\right|\leq\rho^2\left(\frac{2}{h}+\delta\right).
\end{split} 
\end{equation}
For simplicity of notation we can assume that $x=0$ and $\nu=e_2$. With fixed $\eta<1/2$, from (\ref{(16)}) it follows that
\begin{equation}\label{(17)}
\mathcal{A}:=\left|\left((Q_\rho^\nu(x)\cap A_\epsilon)\triangle Q^+_\omega\right)\cap\left\{y:\rho\frac{\eta}{2}\leq\dist(y,\partial Q_\rho^\nu(x))\leq\rho\eta\right\}\right|\leq \rho^2\left(\frac{2}{h}+\delta\right).
\end{equation}
\\
{{\bf Step 3: Construction of an optimal weak circuit.}} We subdivide the annulus between the two concentric squares (with side-lengths $\rho(1-\eta)$ and $\rho(1-2\eta)$ respectively) in four rectangles $R_i(i=1,\dots,4)$ with side-lengths $\rho(1-\eta)$ and $\rho\eta/2$ (they have in common, two by two, a little square of side-length $\rho\eta/2$). From the channel property (Theorem~\ref{channel}), for $\epsilon$ small enough in each of these rectangles we can find at least $c(p)\rho\eta/2\epsilon$ disjoint weak channels $K_i$ connecting the smaller sides of the rectangle, and with length at most $c_1(p)\rho(1-\eta)/\epsilon$. Since $\mathcal{A}\geq\sum_{K_i}(|K_i|\cap\mathcal{A})$, 
from the mean value theorem in each rectangle $R_i$ there exists a weak channel $\widetilde{K}_i$ such that
\begin{equation*}
|\widetilde{K}_i\cap\mathcal{A}|\leq\frac{\mathcal{A}}{\#(K_i)}\leq\frac{\rho^2(2/h+\delta)}{4c(p)\rho\eta/2\epsilon}=\frac{\rho\epsilon}{2c(p)\eta}\left(\frac{2}{h}+\delta\right). 
\end{equation*}
Therefore, considering the weak circuit $K$ contained in $\displaystyle\bigcup_{i=1}^4\widetilde{K}_i$, we have that
\begin{equation}\label{(18)}
\mathcal{H}^1\left(\left((Q_\rho^\nu(x)\cap A_\epsilon)\triangle Q^+_\omega\right)\cap K\right)\leq\frac{2\rho\epsilon}{c(p)\eta}\left(\frac{2}{h}+\delta\right). 
\end{equation}
\\
{{\bf Step 4: Definition of the test sets.}} Now we define the subset $A^1_\epsilon\subset Q_\rho^\nu(x)$ as (see Fig.~\ref{testset})
\begin{equation}
A^1_\epsilon=
\begin{cases}
A_\epsilon,&\text{in the set containing 0 and whose boundary is $K$}\\
Q^+_\omega,&\text{otherwise.}
\end{cases}
\end{equation}

\noindent
Note that
\begin{equation}\label{(20)}
\mathcal{H}^1\left((\partial A^1_\epsilon\backslash\partial A_\epsilon)\cap Q_\rho^\nu(x)\right)\leq\frac{2\rho\epsilon}{c(p)\eta}\left(\frac{2}{h}+\delta\right)+\tilde{c}(p)\rho\delta/2\epsilon. 
\end{equation}
We can find points $x_\epsilon,y_\epsilon\in{\hat{\Z}}^2$ such that $\epsilon x_\epsilon,\epsilon y_\epsilon\in\partial A^1_\epsilon$ and $|\epsilon x_\epsilon+\frac{\rho}{2}e_1|\leq 2\epsilon,|\epsilon y_\epsilon-\frac{\rho}{2}e_1|\leq 2\epsilon$, and a weak path $K^\epsilon$ in $\frac{1}{\epsilon}(\partial A_\epsilon^1\cap Q_\rho^\nu(x))\cap{\hat{\Z}}^2$ connecting $x_\epsilon$ at $y_\epsilon$. By the estimate (\ref{(20)}) we have
\begin{equation*}
\begin{split}
\mu_\epsilon(Q_\rho^\nu(x))&\geq\epsilon|K^\epsilon|-\left(\frac{2\rho\epsilon^2}{c(p)\eta}\left(\frac{2}{h}+\delta\right)+\tilde{c}(p)\rho\delta/2\right)\\
                                                    &\geq\epsilon D^\omega(x_\epsilon,y_\epsilon)-\left(\frac{2\rho\epsilon^2}{c(p)\eta}\left(\frac{2}{h}+\delta\right)+\tilde{c}(p)\rho\delta/2\right).
\end{split}
\end{equation*}
Since $|(y_\epsilon-x_\epsilon)-\frac{\rho}{\epsilon}e_1|\leq 4$, choosing $m=\rho/\epsilon$ in the definition of $\lambda_p$ (equation (\ref{samelim})) and for fixed $\eta, \delta$ and $h$ we obtain
\begin{equation*}
\mathop{\lim\inf}_{\epsilon\to 0}\mu_\epsilon(Q_\rho^\nu(x))\geq\rho\lambda_p(e_1)-\tilde{c}(p)\rho\delta/2=\rho\lambda_p(e_2)-\tilde{c}(p)\rho\delta/2=\rho\lambda_p(\nu)-\tilde{c}(p)\rho\delta/2.
\end{equation*}
By the {\bf(iv)} above we then have
\begin{equation*}
\mathop{\lim\inf}_{\epsilon\to 0}\mu_\epsilon(Q_\rho^\nu(x))\geq\int_{Q_\rho^\nu(x)\cap\partial^*A}\lambda_p(\nu(y))\,d\mathcal{H}^1(y)-\left(\frac{\rho}{h}+\tilde{c}(p)\frac{\rho^2\delta}{2}\right)
\end{equation*}
and finally
\begin{equation*}
\begin{split}
\mathop{\lim\inf}_{\epsilon\to 0}\mu_\epsilon(\Omega)&\geq\sum_j\mathop{\lim\inf}_{\epsilon\to0}\mu_\epsilon(Q_{\rho_j}^{\nu_j}(x_j)\cap\partial^*A)\\
                                                                                    &\geq\sum_j\int_{Q_{\rho_j}^{\nu_j}(x_j)\cap\partial^*A}\lambda_p(\nu(y))\,d\mathcal{H}^1(y)-C\left(\frac{\rho}{h}+\tilde{c}(p)\frac{\rho^2\delta}{2}\right)\\
                                                                                    &=\int_{\Omega\cap\partial^*A}\lambda_p(\nu(y))\,d\mathcal{H}^1(y)-C\left(\frac{\rho}{h}+\tilde{c}(p)\frac{\rho^2\delta}{2}\right);
\end{split}
\end{equation*}
the $\mathop{\lim\inf}$ inequality then follows by the arbitrariness of $\rho,\delta$ and $h$.

The construction of a recovery sequence giving the upper bound (limsup inequality) can be performed just for \emph{polyhedral sets}, since they are dense in energy in the class of sets of finite perimeter. We only give the construction when the set is of the form $\Pi^\nu(x)\cap \Omega$ since it is easily generalized to each face of a polyhedral boundary. {We can localize the construction to the faces of a polyhedral set because the limit energy does not concentrate at its corners: this follows by the chosen scaling.}\\
It is no restriction to suppose that $\Pi^\nu(x)=\Pi^\nu(0)=:\Pi^\nu$, that $\nu$ is a \emph{rational direction} (that is, there exists a positive real number $S$ such that $S\nu\in\Z^2$), and that
\begin{equation}
\mathcal{H}^1(\partial \Omega\cap\partial\Pi^\nu)=0,
\label{(21)}
\end{equation}
since also with these restrictions we obtain a dense class of sets. We will compute the $\Gamma$\hbox{-}limsup for $u=2\chi_{\Pi^\nu}-1$.

Let $M>0$ be large enough so that $\Omega\subset\subset Q_M^\nu(0)$, we set $\tau=\nu^\perp$ and we fix $\eta>0$ such that $\eta<M/2$. There exists a path $\gamma_\epsilon$ in the weak cluster of the dual lattice $\hat{\Z}^2$ contained in the stripe $\{x:|\langle x,\nu\rangle|\leq\eta/\epsilon\}$ and with the two endpoints lying at distance at most $2\epsilon$ from the two sides $\{x:|\langle x,\tau\rangle|=\pm M/2\}$. {The existence of $\gamma_\epsilon$ can be proved with the same construction performed for $\gamma_\omega$ in the proof of the $\Gamma$\hbox{-}$\mathop{\lim\inf}$ inequality.}
After identifying $\gamma_\epsilon$ with a curve in $\R^2$, for $\epsilon$ small enough it disconnects $\frac{1}{\epsilon}\Omega$. We can therefore consider $\Omega_\epsilon^+$, the maximal connected component $\frac{1}{\epsilon}\Omega\backslash\gamma_\epsilon$ containing $\Omega\cup\{\langle x,\nu\rangle\geq\eta/\epsilon\}$, and define
\begin{equation*}
u_\epsilon^\eta(\epsilon i)=
\begin{cases}
1&\text{if $i\in\Z^2\cap \Omega_\epsilon^+$}\\
-1&\text{otherwise.}
\end{cases}
\end{equation*}
Note that
\begin{equation*}
E_\epsilon^\omega(u_\epsilon^\eta)\leq\epsilon|\gamma_\epsilon|\leq\lambda_p(\tau)\mathcal{H}^1(\partial \Pi^\nu\cap \Omega)+O(\eta).
\end{equation*}
\\
By a diagonal argument, for any fixed $\eta>0$ we can construct a subsequence (still denoted by $u_\epsilon^\eta$) converging in $L^1(\Omega)$ as $\epsilon\to0$ to $u^\eta$, where $u^\eta$ is a function such that $\|u^\eta-~u\|_{L^1(\Omega)}\to~0$ as $\eta\to0$.
We have that
\begin{equation*}
\mathop{\lim\sup}_{\epsilon\to0^+} E_{\epsilon}^\omega(u_\epsilon^\eta) \leq \lambda_p(\tau)\mathcal{H}^1(\partial \Pi^\nu\cap \Omega)+O(\eta),\quad\forall\eta>0,
\end{equation*}
and letting $\eta\to0$ we obtain
\begin{equation*}
\Gamma\hbox{-}\mathop{\lim\sup}_{\epsilon\to0^+}E_\epsilon^\omega(u)\leq \lambda_p(\tau)\mathcal{H}^1(\partial \Pi^\nu\cap {\Omega})=\lambda_p(\nu)\mathcal{H}^1(\partial \Pi^\nu\cap {\Omega}).
\end{equation*}
Eventually, we obtain the desired inequality recalling that $\mathcal{H}^1(\partial \Pi^\nu\cap \bar{\Omega})=\mathcal{H}^1(\partial \Pi^\nu\cap \Omega)$ by (\ref{(21)}).

\bigskip
{\bf (b)} It will suffice to show that $\displaystyle E'(u):=\Gamma\hbox{-}\mathop{\lim\inf}_{\epsilon\to0}E_\epsilon^\omega(u)=+\infty$ if $u\not=-1$ or $u\not=1$ identically.  We reason by contradiction and assume that there exists a non-constant function $u\in BV(\Omega;\{\pm1\})$ such that $E'(u)<~+\infty$. Fixed a point $x\in S(u)$ and a square $Q_\rho^\nu(x)$ of side-length $\rho>0$ sufficiently large, by channel property ({see Theorem~\ref{channel} and subsequent remarks}) almost surely there exists (at least) a strong channel connecting two opposite sides of the square. Therefore, if ${u_\epsilon}$ is a sequence converging to $u$, there must be at least one pair $i,j$ of nearest neighbors in the strong cluster such that $(u_\epsilon)_i\neq(u_\epsilon)_j$ so that $E^\omega_\epsilon(u_\epsilon)=+\infty$. This implies that $E'(u)=+\infty$.

\endproof

\section{A continuity result}\label{limitcase}
The number $\lambda_p(\nu)$ defined by equation (\ref{samelim}) describes the average distance on the weak cluster in the direction $\nu$ (and, by Remark~\ref{symmetry}, also in the orthogonal direction). Its value cannot be decreased by using `small portions' of strong connections, as expressed by the following result.
\begin{lem}
Let $\eta>0$ be fixed. Then there are $\delta\in(0,1)$ and $\rho>0$ such that almost surely there exists $N_0$ such that for all $N\geq N_0$ and all channels of length $L$ connecting the two shorter sides of $NT_\nu^\delta$ and with $L<(\lambda_p(\nu)-\eta)N$ we have $\#(\text{strong links})\geq\rho(\eta)N$.
\label{percentage}
\end{lem}
The proof of this technical Lemma is contained in {\scshape Braides} and {\scshape Piatnitski} \cite{BraPia} and it is used to prove that, in the subcritical case, the overall behavior of a discrete membrane with randomly distributed defects is characterized by a fracture type energy, and the surface interaction is described by the asymptotical chemical distance $\lambda_p$. 

We would like to exploit Lemma~\ref{percentage} to prove that an elliptic random spin system with coefficients 1 and $\beta>0$, in the limit as $\beta$ goes to $+\infty$ (i.e. for $\beta$ very large), has the same behavior of a rigid system (that is, with $\beta=+\infty$). More precisely, if $\varphi_p(\beta,\nu)$ is the surface tension coming from the elliptic problem, then

\begin{equation*}
\displaystyle\lim_{\beta\to+\infty}\varphi_p(\beta,\nu)=\lambda_p(\nu)=\varphi_p(+\infty,\nu).
\end{equation*}
The expression $\lambda_p(\nu)=\varphi_p(+\infty,\nu)$ means that $\lambda_p$ is the surface tension computed for $\beta=+\infty$. Such a continuity result seems to be interesting, because in general it does not hold outside this random setting, as shown by the following simple example.

\begin{example}
Consider the energies

\begin{equation}
F_\epsilon^\beta(u)=\sum_{ij}\epsilon c_{ij}^\beta(u_i-u_j)^2,
\end{equation}
where

\begin{equation}
c_{ij}^\beta=
\begin{cases}
\beta & \text{if $i_1=j_1=0$}\\
1      & \text{otherwise.}
\end{cases}
\end{equation}
It is well known from the theory of homogenization of elliptic spin systems that there exists the $\Gamma$\hbox{-}limit

\begin{equation}
\Gamma\hbox{-}\displaystyle\lim_{\epsilon\to0}F_\epsilon^\beta(u)=F^\beta(u)=\displaystyle\int_{{\Omega}\cap\partial^*\{u=1\}}\varphi_\beta(\nu)\,d\mathcal{H}^1=\displaystyle\int_{{\Omega}\cap\partial^*\{u=1\}}\|\nu\|_1\,d\mathcal{H}^1;
\end{equation}
note that the sequences $\{c_{ij}^\beta\}$ and $\{\varphi_\beta\}$ (and consequently $\{F^\beta\}$) are (trivially) increasing in $\beta$ and we can put $c_{ij}^\infty=\displaystyle\sup_{\beta>0}c_{ij}^\beta$, $\widetilde{\varphi}=\displaystyle\sup_{\beta>0}\varphi_\beta$ and

\begin{equation*}
\widetilde{F}(u)=\displaystyle\sup_{\beta>0}F^\beta(u)=\int_{{\Omega}\cap\partial^*\{u=1\}}\widetilde{\varphi}(\nu)\,d\mathcal{H}^1=\displaystyle\int_{{\Omega}\cap\partial^*\{u=1\}}\|\nu\|_1\,d\mathcal{H}^1.
\end{equation*}
\\
Now if we consider the energies

\begin{equation}
F_\epsilon^\infty(u)=\sum_{ij}\epsilon c_{ij}^\infty(u_i-u_j)^2,
\end{equation}
where

\begin{equation}
c_{ij}^\infty=
\begin{cases}
+\infty & \text{if $i_1=j_1=0$}\\
1      & \text{otherwise.}
\end{cases}
\end{equation}
(with the usual convention $+\infty\cdot0=0$), then we have that 

\begin{equation*}
F^\infty(u)=\Gamma\hbox{-}\displaystyle\lim_{\epsilon\to0}F_\epsilon^\infty(u)=\int_{{\Omega}\cap\partial^*\{u=1\}\cap\{{x_1}>0\}}\|\nu\|_1\,d\mathcal{H}^1+\int_{{\Omega}\cap\partial^*\{u=1\}\cap\{{x_1}<0\}}\|\nu\|_1\,d\mathcal{H}^1.
\end{equation*}
\\
Therefore, $\widetilde{F}(u)\neq F^\infty(u)$.
\end{example}

We recall the main results about homogenization of random spin systems (see {\scshape Braides} and {\scshape Piatnitski} \cite{BraPia2}).
Given a probability space $(\Sigma, \mathcal{F}, {\textbf P})$, we consider an ergodic stationary discrete random process $\sigma_{\hat z}^\omega,{\hat z}\in\hat{\Z}^2$.
\\
We are going to compute the $\Gamma$\hbox{-}limit of the energies

\begin{equation*}
E_\varepsilon^\omega(u):=\sum_{ij}\epsilon\sigma_{ij}^\omega (u_i-u_j)^2
\end{equation*}
(with the usual identification $\sigma_{ij}^\omega=\sigma^\omega_{\hat{z}}$). {For any $x,y\in\mathbb{Z}^2$} and $\omega\in\Sigma$ we define

\begin{equation}
\psi^\omega(x,y)=\min\left\{\sum_{n=1}^K\sigma^\omega_{i_n i_{n-1}}:i_0=x,i_K=y,K\in\mathbb{N}\right\},
\label{tension}
\end{equation}
where the minimum is taken over all paths joining $x$ and $y$. The following statement holds (we can compare it with Lemma~\ref{lim}).
\begin{prop}
For any $\tau\in\mathbb{R}^2$ the following limit exists {\bf P}\hbox{-}almost surely and does not depend on $\omega$

\begin{equation}
\psi_0(\tau)=\lim_m\frac{1}{m}\psi^\omega(0,\lfloor m\tau\rfloor),
\label{limtension}
\end{equation}
where $\lfloor m\tau\rfloor_k=\lfloor m\tau_k\rfloor$ is the integer part of the $k$-th component of $m\tau$. \\
For any $x\in\mathbb{R}^2$ and $\tau\in\mathbb{R}^2$ the limit relation

\begin{equation}
\lim_{m\to +\infty}\frac{1}{m}\psi^\omega(\lfloor mx\rfloor,\lfloor mx+m\tau\rfloor)=\psi_0(\tau)
\label{translational}
\end{equation}
holds {\bf P}\hbox{-}almost surely.
\end{prop}

At this point we can recall the main convergence theorem.

\begin{thm}[Elliptic random homogenization]
Let $\sigma_{ij}^\omega$ satisfy the hypothesis of \emph{ellipticity} $0<\alpha\leq \sigma_{ij}^\omega\leq\beta<+\infty$ for all $i,j$. Then the $\Gamma\hbox{-}\displaystyle\lim_{\epsilon\to0}E_\epsilon^\omega$ exists {\bf P}\hbox{-}almost surely, is deterministic and is given by
\begin{equation}
F^\omega(u)=\int_{{\Omega}\cap\partial^*\{u=1\}}\varphi_p(\nu)\,d\mathcal{H}^1.
\label{gammalim}
\end{equation}
where
\begin{equation*}
\label{asint}
\varphi_p(\nu)=\psi_0(\nu^\perp).
\end{equation*}
\label{randomh}
\end{thm}

A particular case of the preceding random problem is obtained by considering the i.i.d. Bernoulli bond\hbox{-}percolation model with coefficients

\begin{equation}
\label{constant}
\sigma^\omega_{\hat{z}}=
\begin{cases}
\beta>0& \mbox{with probability}\quad p\\
1&\mbox{with probability}\quad1-p.
\end{cases}
\end{equation}
With this choice of coefficients we find that the function $\varphi_p$ of Theorem~\ref{randomh} now depends also on $\beta$, $\varphi_p=\varphi_{p,\beta}$. We will prove that the case of rigid spin is the limit as $\beta\to+\infty$ of the problem defined by coefficients (\ref{constant}), in the sense expressed by the following theorem.

\begin{thm}[Continuity] Let $\varphi_{p,\beta}$ be defined by Theorem~\ref{randomh} when the coefficients $\sigma^\omega_{\hat{z}}$ are given by (\ref{constant}), and $\lambda_p$ be defined by Proposition~\ref{generallim}. Then, for all $\nu\in\R^2$, we have two cases:
\begin{description}
\item[(i)] If $p<1/2$, then $\displaystyle\lim_{\beta\to+\infty}\varphi_{p,\beta}(\nu)=\lambda_p(\nu)$;
\item[(ii)] If $p>1/2$, then $\displaystyle\lim_{\beta\to+\infty}\varphi_{p,\beta}(\nu)=+\infty$.
\end{description}
\label{continuo}
\end{thm}
\proof
{\bf (i)} Let $p<1/2$. First remark that, with fixed $\nu\in\R^2$ and for all $\beta>0$, we have
\begin{equation*}
\varphi_{p,\beta}(\nu)\leq\lambda_p(\nu),
\end{equation*}
because the minimum in the definition of $\lambda_p$ (see equation (\ref{samelim})) is taken in a smaller set of paths.\\
Therefore, $\lambda_p$ being independent of $\beta$,
\begin{equation*}
\displaystyle\lim_{\beta\to+\infty}\varphi_{p,\beta}(\nu)\leq\lambda_p(\nu).
\end{equation*}
Now suppose that there exists $\eta>0$ such that $\displaystyle\lim_{\beta\to+\infty}\varphi_{p,\beta}(\nu)\leq\lambda_p(\nu)-\eta$; we would like to show that this assumption leads to a contradiction. 
\\
First of all, $\varphi_{p,\beta}$ being increasing in $\beta$, we have that
\begin{equation}
\varphi_{p,\beta}(\nu)\leq\lambda_p(\nu)-\eta,\quad \forall\beta>0.
\end{equation}
With fixed $\beta$, by (\ref{translational}) there exists $\bar{n}\in\N$ almost surely such that, for all $n\geq\bar{n}$ 
\begin{equation}
\psi^\omega(\lfloor nx\rfloor,\lfloor nx+n\tau\rfloor)\leq(\lambda_p(\nu)-\eta')n,
\label{(4.8)}
\end{equation}
where $\tau=\nu^\perp$, $\eta'$ is a constant and $x\in\R^2$.
Suppose that $\displaystyle\psi^\omega(\lfloor nx\rfloor,\lfloor nx+n\tau\rfloor)=\sum_{m=1}^K\widetilde{\sigma}^\omega_{i_m i_{m-1}}$ with $ i_0=\lfloor nx\rfloor,i_K=\lfloor nx+n\tau\rfloor$ and let $\gamma$ be the corresponding path.\\
{By means of Lemma~\ref{percentage}, we can find $\delta\in(0,1)$, a constant $C=C(\eta')>0$ such that if $n'\geq n$ is such that the channel $\gamma$ connects the shorter sides of $n'T^\delta_\nu$, then from the fact that
\begin{equation}
|\gamma|\leq\sum_{m=1}^K\widetilde{\sigma}^\omega_{i_m i_{m-1}}\leq(\lambda_p(\nu)-\eta')n',
\end{equation}}
\\
{it follows that
\begin{equation}
\#(\text{strong links in }\gamma)\geq Cn'.
\end{equation}}
{Now
\begin{equation}
\beta Cn'\leq\sum_{m=0}^K\widetilde{\sigma}^\omega_{i_m i_{m-1}}\leq(\lambda_p(\nu)-\eta')n',
\end{equation}
and letting $\beta\to+\infty$ we obtain a contradiction.}

\bigskip
{\bf (ii)} If $p>1/2$, we can argue as in Theorem~\ref{percthm}{\bf (b)}, because the percentage of $\beta$ is fixed by channel property. In particular, for large $m$, the paths linking $\lfloor mx\rfloor$ and $\lfloor mx+m\tau \rfloor$ in equations (\ref{tension}), (\ref{limtension}) and (\ref{translational}) contain at least a $\beta$\hbox{-}bond.
\endproof

\bigskip
\noindent{\bf Acknowledgements.} I am very grateful to Andrea Braides for suggesting this problem, and I would like to thank him for his fundamental advices. I acknowledge Andrey Piatnitski for his useful suggestions for the proof of Proposition~\ref{generallim}, and the anonymous referees for their interesting remarks leading to improvements of the manuscript.

\end{document}